\documentclass[11pt]{article}

\newenvironment{proof}{\noindent {\em Proof}.\ }{\hspace*{\fill}$\halmos$\medskip}
\newcommand{\halmos}{\rule{1ex}{1.4ex}}

\newtheorem{theorem}{Theorem}

\newtheorem{remark}{Remark}

\newcommand{\enma}[1]   {\ensuremath{#1}}

\newcommand{\beq}{\begin{equation}}
\newcommand{\eeq}{\end{equation}}
\newcommand{\beqn}{\begin{eqnarray}}
\newcommand{\eeqn}{\end{eqnarray}}
\newcommand{\bseq}{\begin{subequations}}
\newcommand{\eseq}{\end{subequations}}
\newcommand{\ba}{\begin{array}}
\newcommand{\ea}{\end{array}}
\newcommand{\bct}{\begin{center}}
\newcommand{\ect}{\end{center}}
\newcommand{\btmz}{\begin{itemize}}
\newcommand{\etmz}{\end{itemize}}
\newcommand{\benum}{\begin{enumerate}}
\newcommand{\eenum}{\end{enumerate}}

\newcommand{\diag}      {\enma{\mathrm{diag}}}

\newcommand{\sgn}       {\enma{\mathrm{sgn}}}

\newcommand{\inner}[2]{\left\langle #1,#2 \right\rangle}

\newcommand{\matbegin}{
        \left[
}
\newcommand{\matend}{
        \right]
}

\newcommand{\obth}[3]{
  \matbegin \begin{array}{ccc}
       #1 & #2 & #3
       \end{array} \matend }

\newcommand{\ca}{{\cal A}}

\newcommand{\be}{\begin{equation}}
\newcommand{\ee}{\end{equation}}

\newcommand{\cplxs}{ C\kern -.35em \rule{0.03 em}{.7 ex}~   }

\def\complex{\hbox{C\kern -.45em \rule{0.03 em}{1.5 ex}}~}

\newcommand{\bi}{\begin{itemize}}
\newcommand{\ei}{\end{itemize}}

\newtheorem{assumption}{Assumption}

\usepackage{amssymb,latexsym,color,amsmath,pifont,epsfig,cite}
\usepackage{setspace,mathtools}

     \mathtoolsset{showmanualtags}

\newcommand{\bbN}{\mathbb{N}}
\newcommand{\bbR}{\mathbb{R}}

\newcommand{\non}{\nonumber}
\newcommand{\ds}{\displaystyle}

\newcommand{\mrd}{\mathrm{d}}
\newcommand{\mre}{\mathrm{e}}

\newcommand{\Dom}{\Omega}

\begin{document}

\title{\bf \LARGE A passivity-based approach to stability of spatially distributed
systems with a cyclic interconnection structure}

\author{Mihailo R.\ Jovanovi\'c, Murat Arcak, and Eduardo D.\ Sontag
\thanks{M.\ R.\ Jovanovi\'c is with the Department of Electrical
and Computer Engineering, University of
Minnesota,~Minneapolis,~MN~55455, USA (mihailo@umn.edu).}
\thanks{M.\ Arcak is with the Department of Electrical, Computer,
and Systems Engineering, Rensselaer Polytechnic
Institute,~Troy,~NY~12180, USA \mbox{(arcakm@rpi.edu)}.}
\thanks{E.\ D.\ Sontag is with the Department of Mathematics,
Rutgers University, Piscataway, NJ~08854, USA
(sontag@math.rutgers.edu).}
}%
    \maketitle

    \begin{abstract}
A class of distributed systems with a cyclic interconnection
structure is considered. These systems arise in several biochemical
applications and they can undergo diffusion driven instability which
leads to a formation of spatially heterogeneous patterns. In this
paper, a class of cyclic systems in which addition of diffusion does
not have a destabilizing effect is identified. For these systems
global stability results hold if the ``secant'' criterion is
satisfied. In the linear case, it is shown that the secant condition
is necessary and sufficient for the existence of a decoupled
quadratic Lyapunov function, which extends a recent diagonal
stability result to partial differential equations. For
reaction-diffusion equations with nondecreasing coupling
nonlinearities global asymptotic stability of the origin is
established. All of the derived results remain true for both linear
and nonlinear positive diffusion terms. Similar results are shown
for compartmental systems.
    \end{abstract}


\section{Introduction}

The first gene regulation system to be studied in detail was the one
responsible for the control of lactose metabolism in \emph{E. Coli},
the \emph{lac} operon studied in the classical work of Jacob and Monod
\cite{JacobMonod61,MonodChangeauxJacob63}.
Jacob and Monod's work led Goodwin \cite{Goodwin65} and later many others
\cite{%
FraserTiwari74, TiwariFraserBeckmann74, AaronsGray75,
HasTysWeb77, SanglierNicolis76, tysoth78, BanksMahaffy781,
BanksMahaffy782, Mahaffy84, Mahaffy842, Glass-oscill, thr91}
to the mathematical study of systems made up of cyclically interconnected
genes and gene products.
In addition to gene regulation networks, cyclic feedback structures have been
used as models of certain metabolic pathways \cite{MorKay67}, of tissue growth
regulation \cite{WeiKav57},
of cellular signaling pathways \cite{KholodenkoNegfeedback},
and of neuron models \cite{Hei76}.

Generally, cyclic feedback systems (of arbitrary order) were shown by
Mallet-Paret and Smith \cite{malletparet-smith, Mallet-ParetDelay} to have
behaviors no more complicated that those of second-order systems: for
precompact trajectories, $\omega$-limit sets can only consist of equilibria,
limit cycles, or heteroclinic or homoclinic connections, just as in the planar
Poincar\'e-Bendixson Theorem.
When the net effect around the loop is positive, no (stable)
oscillations are possible,
because the overall system is monotone \cite{smith}.
On the other hand, inhibitory or ``negative feedback'' loops
give rise to the possibility of periodic orbits, and it is then of
interest to provide conditions for oscillations or lack thereof.

Besides the scientific and mathematical interest of the study of cyclic
negative feedback systems, there is an engineering motivation as well, which
rises from the field of synthetic biology.  Oscillators will be fundamental
parts of engineered gene bacterial networks, used to provide timing and
periodic signals to other components.  A major experimental effort,
pioneered by the construction of the ``repressilator'' by Elowitz and Leibler
\cite{Elowitz}, is now under way to build reliable oscillators with gene
products.
Indeed, the theory of cyclic feedback systems has been proposed as a way
to analyze the repressilator and similar systems
\cite{ddv-hesACC,ACC2007-sub}.

In order to evaluate stability properties of negative feedback cyclic systems,
\cite{tysoth78} and \cite{thr91} analyzed the Jacobian linearization at the
equilibrium, which is of the form

    \begin{equation}
    A \; = \;
    \left[\begin{array}{ccccc} -a_1  & 0 & \cdots & 0 & -b_n \\ b_1 & -a_2 &\ddots & & 0 \\
    0 & b_2 & -a_3 & \ddots & \vdots \\ \vdots & \ddots & \ddots &
    \ddots & 0 \\ 0 & \cdots & 0 & b_{n-1} & -a_n
    \end{array} \right]
    \label{Amatrix}
    \end{equation}
$a_i>0,
    ~
    b_i>0,
    ~
    i=1,\cdots,n$, and showed that $A$ is Hurwitz if the following
    sufficient condition
holds:
    \begin{equation}
    \label{eq.secant}
    \frac{b_1\cdots b_n}{a_1\cdots a_n} \, < \, \sec(\pi/n)^n.
    \end{equation}
This ``secant criterion''is also necessary for stability when the
$a_i$'s are identical.

An application of the secant condition in a ``systems biology''
context was in Kholodenko's \cite{KholodenkoNegfeedback} (see
also~\cite{stas2}) analysis of a simplified model of negative
feedback around MAPK (mitogen activated protein kinase) cascades.
MAPK cascades constitute a highly conserved eukaryotic pathway,
responsible for some of the most fundamental processes of life such
as cell proliferation and growth
\cite{ferrell,lauffenburger,widman}. Kholodenko used the secant
condition to establish conditions for asymptotic stability.

\subsection{Global stability considerations}

\begin{figure}[ht]
\centering
\includegraphics[scale=0.35]{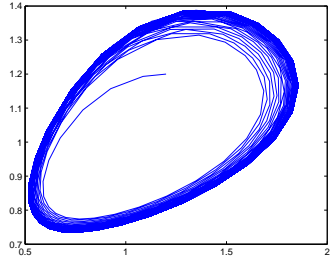}
\caption{Trajectory of (\ref{counterex}) starting from initial
condition $\chi=[1.2\ 1.2\ 1.2]^T$, projected onto the
$\chi_1$-$\chi_2$ plane.} \label{USE}
\end{figure}

It appears not to be generally appreciated that (local) stability of the
equilibrium in a cyclic negative feedback system does not rule out
the possibility of periodic orbits. Indeed, the
Poincar\'{e}-Bendixson Theorem of Mallet-Paret and Smith
\cite{malletparet-smith, Mallet-ParetDelay} allows such periodic
orbits to coexist with stable equilibria. As an illustration
consider the system
\begin{eqnarray}\nonumber
\dot{\chi}_1&=&-\chi_1+\varphi(\chi_3)\\
\dot{\chi}_2&=&-\chi_2+\chi_1 \label{counterex} \\
\dot{\chi}_3&=&-\chi_3+\chi_2 \nonumber
\end{eqnarray}
where
\begin{equation}\label{sat}
\varphi(\chi_3)=e^{-10(\chi_3-1)}+0.1{\rm sat}(25(\chi_3-1)),
\end{equation}
and ${\rm sat}(\cdot):=\sgn(\cdot)\min\{1,|\cdot|\}$ is a
saturation\footnote{One can easily modify this example to make
$\varphi(\cdot)$ smooth while retaining the same stability
properties.} function. The function (\ref{sat}) is decreasing, and
its slope has magnitude $b_3=7.5$ at the equilibrium
$\chi_1=\chi_2=\chi_3=1$. With $a_1=a_2=a_3=b_1=b_2=1$ and $n=3$,
the secant criterion (\ref{eq.secant}) is satisfied and, thus, the
equilibrium is asymptotically stable. However, simulations in
Fig.~\ref{USE} show the existence of a periodic orbit in addition to
this stable equilibrium.

To delineate {\em global} stability properties of cyclic systems
with negative feedback, \cite{arcson06} studied (by building on a
passivity interpretation of the secant criterion in \cite{secant})
the nonlinear model
    \beq
    \ba{rcl}
    \dot{x}_1
    & \!\! = \!\! &
    -f_1(x_1) \, - \, g_n(x_n)
    \\[0.1cm]
    \dot{x}_2
    & \!\! = \!\! &
    -f_2(x_2) \, + \, g_1(x_1)
    \\[0.1cm]
    &\vdots&
    \\[0.1cm]
    \dot{x}_n
    & \!\! = \!\! &
    -f_n(x_n) \, + \, g_{n-1}(x_{n-1})
    \label{bio}
    \ea
    \eeq
and proved global asymptotic stability of the origin\footnote{In the
rest of the paper we assume that an equilibrium exists and is unique
(see \cite{arcson06} for conditions that guarantee this) and that
this equilibrium has been shifted to the origin with a change of
variables.} under the conditions
    \begin{equation}
    \sigma f_i(\sigma) \, > \, 0,
    ~~
    \sigma g_i(\sigma) \, > \, 0,
    ~~
    \forall \,
    \sigma
    \, \in \,
    \bbR \setminus \{0\},
    \tag*{(C1)}
    \label{eq.C1}
    \end{equation}
    \begin{equation}
    \frac{g_{i}(\sigma)}{f_i(\sigma)}
    \, \leq \,
    \gamma_i,
    ~~
    \forall \,
    \sigma
    \, \in \,
    \bbR \setminus \{0\},
    \tag*{(C2)}
    \label{eq.C2}
    \end{equation}
    \begin{equation}
    {\gamma_1 \cdots \gamma_n}
    \, < \,
    \sec(\pi/n)^n,
    \tag*{(C3)}
    \label{eq.C3}
    \end{equation}
    \begin{equation}
    \lim_{|x_i| \, \rightarrow \, \infty}
    \int_0^{x_i} g_i(\sigma) \, \mrd
    \sigma
    \, = \, \infty.
    \tag*{(C4)}
    \label{eq.C4}
    \end{equation}
The conditions \ref{eq.C1}-\ref{eq.C4} encompass the linear system
(\refeq{Amatrix})-(\refeq{eq.secant}) in which $f_i(x_i)=a_ix_i$,
$g_i(x_i)=b_ix_i$, and $\gamma_i=b_i/a_i$.

A crucial ingredient in the global asymptotic stability
 proof of \cite{arcson06} is the observation that the secant
 criterion (\ref{eq.secant}) is necessary and sufficient for
 {\em diagonal stability\/} of (\ref{Amatrix}), that is for the existence of
 a diagonal matrix $D>0$ such that
 \begin{equation}
 \label{dstable}
 A^T D \; + \; D A \; < \; 0.
 \end{equation}
Using this diagonal stability property, \cite{arcson06} constructs a
Lyapunov function for (\ref{bio}) which consists of a weighted sum
of decoupled functions of the form $V_i(x_i)=\int_0^{x_i}g_i(\sigma)
\, \mrd \sigma$. In the linear case this construction coincides with
the quadratic Lyapunov function $V=x^TDx$.

\subsection{Spatial localization}

Ordinary differential equation models such as described above implicitly
assume that reactions proceed in a ``well-mixed'' environment.
However, in cells, certain processes are localized to membranes (activation of
pathways by receptors), to the nucleus (transcription factor binding to DNA,
production of mRNA), to the cytoplasm (much of signaling), or to one of the
specialized organelles in eukaryotes.
The exchange of chemical species between these spatial domains has been found
to be responsible for dynamical behavior, such as emergence of oscillations,
in fundamental cell signaling pathways, see for instance \cite{hoffmann02}.
These exchanges often happen by random movement (diffusion), although
transport mechanisms and gated channels are sometimes involved as well.

When each of a finite set of spatial domains is reasonably
``well-mixed,'' so that the concentrations of relevant chemicals in
each domain are appropriately described by ordinary differential
equations (ODEs), a compartmental model may be used.  In a
compartmental model, several copies of an ODE system are
interconnected by ``pipes'' that tend to balance species
concentrations among connected compartments.  The overall system is
still described by a system of ODEs, but new dynamical properties
may emerge from this interconnection. For example, two copies of an
oscillating system may synchronize, or two multi-stable systems may
converge to the same steady state.

On the other hand, if a well-mixed assumption in each of a finite
number of compartments is not reasonable, a more appropriate
mathematical formalism is that of reaction-diffusion partial
differential equations (PDEs)~\cite{amann1,hen81,pazy,rot84,smo94}:
instead of a dynamics $\dot x=f(x)$, one considers equations of the
general form \beq \frac{\partial x}{\partial t} = D\Delta x +
f(x)\,,
\quad\quad
\frac{\partial x}{\partial \nu } = 0\,,
\label{eq:pde1} \eeq where now the vector $x=x(\xi,t)$ depends on
both time $t$ and space variables $\xi$ belonging to some domain
$\Dom$, $\Delta x$ is the Laplacian of the vector $x$ with respect
to the space variables, $D$ is a matrix of positive diffusion
constants, and ${\partial x}/{\partial \nu }$ denotes the
directional derivative in the direction of the normal to the
boundary $\partial \Dom$ of the domain $\Dom$, representing a
no-flux or Neumann boundary condition. (Technical details are given
later, including generalizations to more general elliptic operators
that model space-dependent diffusions.)

Diffusion plays a role in generating new behaviors for the PDE as compared to
the original ODE $\dot x=f(x)$.
In fact, one of the main areas of research in mathematical biology concerns
the phenomenon of diffusive instability, which constitutes the basis of
Turing's mechanism for pattern
formation~\cite{turing,othmer-turing,murray89,keshet05},
and which amounts to the emergence of stable non-homogeneous in space
solutions of a reaction-diffusion PDE.
The Turing phenomenon has a simple analog, and is easiest to
understand intuitively, for an ODE consisting of two identical compartments
\cite{keshet05,iglesias_compartment}.
Also in the context of cell signaling, and in particular for the MAPK pathway
mentioned earlier, reaction-diffusion PDE models play an important
role~\cite{kholodenkodiffusion06}.

If diffusion coefficients are very large, diffusion effects may be ignored in
modeling.
As an illustration, the stability of uniform steady states is unchanged
provided that the diffusion coefficient $D$ is sufficiently large compared to
the ``steepness'' of the reaction term $f$, measured for instance by an upper
bound $a$ on its Lipschitz constant or equivalently the maximum of its
Jacobians at all points (for chemical reaction networks, this is interpreted
as the inverse of the kinetic relaxation time, for steady states).
Introducing an energy function using the integral of
$|\partial x/\partial  \xi |^2$, and then integrating by parts and using
Poincar\'e's inequality,
one obtains an exponential decrease of this energy, controlled by
the difference of $a$ and $D$ (\cite{jones-sleeman}, Chapter 11).
For instance, Othmer~\cite{othmer-large-diffusion} provides a
condition $D\mu >a$ in terms of the smallest nonzero eigenvalue of
the Neumann Laplacian $\{ \Delta x + \mu x=0, \,\xi \in \Dom$;
$\partial x/\partial \nu = 0, \,\xi \in  \partial \Dom \}$ to
guarantee exponential convergence to zero of spatial
nonuniformities, and estimates that his condition is met for
intervals $\Omega =[0,L]$ of length $L\approx10\mu m$, with
diffusion of at least about 4$\times 10^{-8}$cm$^2$/sec and
$a\approx10^{-1}$sec.

\begin{figure}[ht]
    \begin{center}
    \includegraphics[height=2.05in,width=3.in]{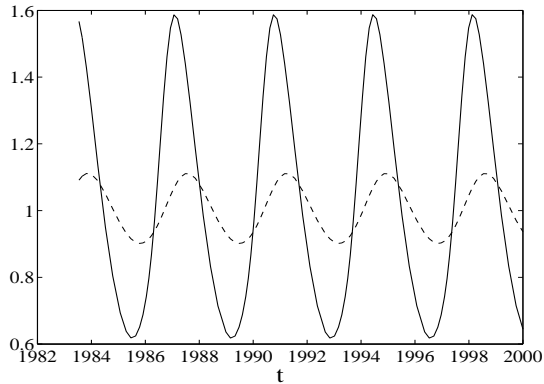}
    \end{center}
    \caption{New oscillations in two-compartment system: $\chi_1$
    (solid) and $\eta_1$ (dashed) shown.}
    \label{fig:2compartmentcounter}
    \end{figure}

On the other hand, if diffusion is not dominant, it is necessary to explicitly
incorporate spatial inhomogeneity, whether through compartmental or PDE models.
The goal of this paper is to extend the linear and nonlinear secant condition
to such compartmental and PDE models, using a passivity-based approach.
To illustrate why spatial behavior may lead to interesting new phenomena even
for cyclic negative feedback systems, we take a two-compartment version
of the system shown in (\ref{counterex}):
\begin{eqnarray}\nonumber
\dot{\chi}_1&=&-\chi_1+\varphi(\chi_3)  + D(\eta_1-\chi_1)\\
\dot{\chi}_2&=&-\chi_2+\chi_1 \nonumber + D(\eta_2-\chi_2)\\
\dot{\chi}_3&=&-\chi_3+\chi_2 \nonumber + D(\eta_3-\chi_3)\\
\dot{\eta}_1&=&-\eta_1+\varphi(\eta_3)  + D(\chi_1-\eta_1)\\
\dot{\eta}_2&=&-\eta_2+\eta_1 \nonumber + D(\chi_2-\eta_2)\\
\dot{\eta}_3&=&-\eta_3+\eta_2 \nonumber + D(\chi_3-\eta_3)
\end{eqnarray}
and pick $D=10^{-4}$. We simulated this system with initial
condition $\left[\,1.4945~1.3844~1.0877~1~1~1\,\right]^T$, so that
the first-compartment $\chi_i(0)$ coordinates start approximately on
the limit cycle, and the second-compartment $\eta_i(0)$ coordinates
start at the equilibrium.  The resulting simulation shows that a new
oscillation appears, in which both components oscillate, out of
phase (no synchronization), with roughly equal period but very
different amplitudes.  Figure~\ref{fig:2compartmentcounter} shows
the solution coordinates $\chi_1$ and $\eta_1$
plotted on a window after a transient behavior.
This oscillation is an emergent behavior of the compartmental system, and is
different from the limit cycle in the original three-dimensional system.
(One may analyze the existence and stability of these orbits using an ISS-like
small-gain theorem.)

Our goal is to show that--in contrast to this example--if the secant
condition {\em does\/} apply to a negative cyclic feedback system,
then no non-homogeneous limit behavior can arise, in compartmental
or in PDE models, no matter what is the magnitude of the diffusion
effect.

\section{Problem formulation}

In this paper we extend the linear and nonlinear results of
\cite{tysoth78,thr91,arcson06} to spatially distributed models that
consist of a cyclic interconnection of $n$ reaction-diffusion
equations
    \beq
    \ba{rcl}
    \!\!\!\!\!\!
    \psi_{1t}
    & \!\! = \!\! &
    \nabla \cdot \left( h_1(\psi_1) \, \nabla \psi_{1} \right)
    \, - \,
    f_1 (\psi_1)
    \, - \,
    g_n (\psi_n)
    \\[0.1cm]
    \!\!\!\!\!\!
    \psi_{2t}
    & \!\! = \!\! &
    \nabla \cdot \left( h_2 (\psi_2) \, \nabla \psi_{2} \right)
    \, - \,
    f_2 (\psi_2)
    \, + \,
    g_1 (\psi_1)
    \\[0.1cm]
    \!\!\!\!\!\!
    & \!\! \vdots \!\! &
    \\[0.1cm]
    \!\!\!\!\!\!
    \psi_{nt}
    & \!\! = \!\! &
    \nabla \cdot \left( h_n (\psi_n) \, \nabla \psi_{n} \right)
    \, - \,
    f_n (\psi_n)
    \, + \,
    g_{n-1} (\psi_{n-1})
    \ea
    \tag*{(RD)}
    \label{eq.RD}
    \eeq
where $\psi_i$ denotes the state of the $i$th subsystem which
depends on spatial coordinate $\xi$ and time $t$, $\psi_i(\xi,t)$,
and $f_i$, $g_i$, $h_i$ denote static nonlinear functions of their
arguments. We consider a situation in which the spatial coordinate
$\xi := (\xi_1, \ldots, \xi_r)$ belongs to a bounded domain $\Omega$
in $\bbR^r$, $r = 1$, $2$ or $3$, with a smooth boundary $\partial
\Omega$ and outward unit normal $\nu$. The state of each subsystem
satisfies the Neumann boundary conditions, $\partial
\psi_{i}/\partial \nu := \psi_{i \nu} = 0$ on $\partial \Omega$,
$\nabla \psi_i$ is the gradient of $\psi_i$, $\nabla \cdot v$ is the
divergence of a vector $v$, and the domain of the
\mbox{$r$-dimensional} Laplacian $\Delta := \nabla \cdot \nabla$ is
given by~\cite{ban83,curzwa95}
    \beq
    {\cal D} (\Delta)
    \, := \,
    \left\{
    \psi_i \in H_2 (\Omega),
    ~
    \psi_{i \nu}
    \, = \, 0
    ~\mbox{on}~
    \partial \Omega
    \right\}.
    \tag*{(DM)}
    \label{eq.DM}
    \eeq
Here, $H_2 (\Omega)$ denotes a Sobolev space of square integrable
functions with square integrable second distributional derivatives.
The standard $L_2^n (\Omega)$ inner product is given by
    \beq
    \inner{\psi}{\phi}
    \; := \;
    \int_{\Omega}
    \psi^T (\xi) \, \phi (\xi)
    \,
    \mrd \xi
    \non
    \eeq
where $\mrd \xi \, := \, \mrd \xi_1 \cdots \mrd \xi_r$ and $\psi \,
:= \, \obth{\psi_1}{\cdots}{\psi_n}^T$.

As explained in the introduction, the study of stability properties
for distributed system \ref{eq.RD} is important in many biological
applications. Our first result, presented in Section
\ref{sec.linear}, studies the linearization of \ref{eq.RD} and shows
that the secant condition (\ref{eq.secant}) is sufficient for the
exponential stability despite the presence of diffusion terms. It
further shows that the secant condition is necessary and sufficient
for the existence of a decoupled Lyapunov function, thus extending
the diagonal stability result of \cite{arcson06} to partial
differential equations. The next result of the paper, presented in
Section \ref{NL}, studies the nonlinear reaction-diffusion equation
\ref{eq.RD} and proves global asymptotic stability of $\psi=0$ under
assumptions that mimic the conditions \ref{eq.C1}-\ref{eq.C3} of
\cite{arcson06}, and under the additional assumptions that the
functions $g_i(\cdot)$ and $h_i(\cdot)$, $i=1,\cdots,n$, be
nondecreasing and positive, respectively. This additional assumption
on the $g$-functions ensures convexity of the Lyapunov function
which is a crucial property for our stability proof. Indeed, a
similar convexity assumption has been employed in~\cite{FitHolMor97}
to preserve stability in the presence of linear diffusion terms.
Finally, Section \ref{sec.compart} studies a compartmental ordinary
differential equation model instead of the partial differential
equation \ref{eq.RD}, and proves global asymptotic stability using
the same nondecreasing assumption for $g_i$'s.

\section{Cyclic interconnection of linear reaction-diffusion equations}
    \label{sec.linear}

We start our analysis by considering an interconnection of spatially
distributed systems~\refeq{eq.RD} with
    \beq
    f_i (\psi_i)
    \, := \,
    a_i \psi_i,
    ~
    g_i (\psi_i)
    \, := \,
    b_i \psi_i,
    ~
    h_i (\psi_i)
    \, := \,
    c_i,
    ~
    i \, = \, 1, \ldots, n
    \label{eq.lin-constr}
    \eeq
where each $a_i$, $b_i$, and $c_i$ represents a positive parameter.
In this case, system~\refeq{eq.RD} simplifies to a cascade connection
of linear reaction-diffusion equations where the output of the last
subsystem is brought to the input of the first subsystem through a
negative unity feedback. Abstractly, the dynamics of
system~\refeq{eq.RD}-\refeq{eq.DM} with $f_i(\cdot)$, $g_i(\cdot)$,
and $h_i(\cdot)$ satisfying~(\refeq{eq.lin-constr}) are given by
    \beq
    \psi_t
    \; = \;
    {\cal A} \psi
    \; := \;
    C \Delta \psi
    \, + \,
    A_0 \psi
    \tag*{(LRD)}
    \label{eq.LRD}
    \eeq
where $\Delta \psi$ denotes the vector Laplacian, that is $\Delta
\psi \, := \, \obth{\Delta \psi_1}{\cdots}{\Delta \psi_n}^T$, $C \,
:= \, \diag \{ \obth{c_1}{\cdots}{c_n} \} \, > \, 0$, and
    \beq
    \ba{c}
    A_0
    \; := \;
    \left[
    \ba{ccccc}
    - \, a_1 & 0 & \cdots & 0 & - \, b_n
    \\
    b_1 & - \, a_2 & \ddots & & 0
    \\
    0 & b_2 & - \, a_3 & \ddots & \vdots
    \\
    \vdots & \ddots & \ddots & \ddots & 0
    \\
    0 & \cdots & 0 & b_{n-1} & - \, a_n
    \ea
    \right],
    ~~
    a_i \, > \, 0,
    ~
    b_i \, > \, 0,
    ~
    i \, = \, 1, \ldots, n.
    \ea
    \non
    \eeq

\subsection{Exponential stability and the secant criterion in one spatial dimension}

In this section, we focus on systems with one spatial dimension $\xi
\in \Omega := (0,1)$. We show that operator $\cal A$
with~\refeq{eq.DM} generates an exponentially stable strongly
continuous ($C_o$) semigroup $T(t)$ on $L_2^n (0,1)$ if the {\em
secant criterion\/} (\refeq{eq.secant}) is satisfied. We note that
the exponential stability of $T(t)$ in
Theorem~\ref{th.heat-exp-stab} can be also established using a
Lyapunov based approach that we develop for systems with two or
three spatial coordinates. However, the proof of
Theorem~\ref{th.heat-exp-stab} is of independent interest because of
the explicit construction of the $C_o$-semigroup and
block-diagonalization of operator~\refeq{eq.LRD}-\refeq{eq.DM}
(which is well suited for a modal interpretation of stability
results in one spatial coordinate).

It is well known (see, for example~\cite{curzwa95}) that the
operator $\partial_{\xi \xi}$ with Neumann boundary conditions is
self-adjoint with the following set of eigenfunctions
$\{\varphi_k\}$ and corresponding eigenvalues $\{\nu_k\}$:
    \[
    \ba{rclrcll}
    \varphi_0(\xi)
    & \!\! = \!\! &
    1,
    &
    \varphi_l(\xi)
    & \!\! = \!\! &
    \sqrt{2 }\cos \, l \pi \xi,
    &
    l \, \in \, \bbN,
    \\[0.1cm]
    \nu_0
    & \!\! = \!\! &
    0,
    &
    \nu_l
    & \!\! = \!\! &
    - (l \pi)^2,
    &
    l \, \in \, \bbN.
    \ea
    \]
Since the eigenfunctions $\{\varphi_k\}$ represent an orthonormal
basis of $L_2 (0,1)$ each $\psi_i(\xi,t)$ can be represented as
    \beq
    \psi_i(\xi,t)
    \; = \;
    \sum_{k \, = \, 0}^{\infty}
    x_{i,k} (t) \varphi_k (\xi),
    \non
    \eeq
where $x_{i,k} (t)$ denote the spectral coefficients given by
    \beq
    x_{i,k} (t)
    \; = \;
    \inner{\varphi_k}{\psi_i}
    \; := \;
    \int_{0}^{1}
    \varphi_k (\xi) \psi_i (\xi,t) \, \mrd \xi.
    \non
    \eeq
Thus, a spectral decomposition of operator $\partial_{\xi \xi}$
in~\refeq{eq.LRD} yields the following infinite-dimensional system
on $l_2^n$ of decoupled $n$th order equations:
    \beq
    \dot{x}_k
    \; = \;
    A_k x_k,
    ~~
    k \, = \, 0, 1, \ldots,
    \label{eq.xk-dot}
    \eeq
with $x_k (t) \, := \, \obth{x_{1,k}(t)}{\cdots}{x_{n,k} (t)}^T$,
    \beq
    A_k
    \, := \,
    \left[
    \ba{ccccc}
    - \, \alpha_{1,k} & 0 & \cdots & 0 & - \, b_n
    \\
    b_1 & - \, \alpha_{2,k} & \ddots & & 0
    \\
    0 & b_2 & - \, \alpha_{3,k} & \ddots & \vdots
    \\
    \vdots & \ddots & \ddots & \ddots & 0
    \\
    0 & \cdots & 0 & b_{n-1} & - \, \alpha_{n,k}
    \ea
    \right],
    \non
    \eeq
and
    $
    \alpha_{i,k}
    \, := \,
    a_i \, - \, c_i \nu_k
    \, = \,
    a_i \, + \, c_i (k \pi)^2
    \, > \,
    0.
    $
Based on~\cite{tysoth78,thr91} we conclude that each $A_k$ is
Hurwitz if~(\refeq{eq.secant}) holds. Therefore, each subsystem
in~(\refeq{eq.xk-dot}) is exponentially stable and there exist $P_k
= P_k^T > 0$ such that
    \beq
    A_k^T P_k
    \, + \,
    P_k A_k
    \, = \,
    - I,
    ~~
    k \, = \, 0, 1, \ldots.
    \non
    \eeq
Now, since $\ca$ is the infinitesimal generator of the following
$C_o$-semigroup:
    \beq
    T(t) \psi (0)
    \; := \;
    T(t) \psi (\xi,0)
    \; = \;
    \sum_{k \, = \, 0}^{\infty}
    \mre^{A_k t} x_k (0) \varphi_k (\xi),
    \non
    \eeq
we have
    \beq
    \ba{rcl}
    \ds{\int_{0}^{\infty}}
    \| T(t) \psi (0) \|^2
    \,
    \mrd t
    \; &:=& \;
    \ds{\int_{0}^{\infty}}
    \inner{T(t) \psi (0)}{T(t) \psi (0)}
    \,
    \mrd t
    \; = \;
    \ds{\sum_{k \, = \, 0}^{\infty}}
    x_k^T (0)
    \left(
    \int_{0}^{\infty}
    \mre^{A_k^T t} \mre^{A_k t}
    \,
    \mrd t
    \right)
    x_k (0)\\
    \; &=& \;
    \ds{\sum_{k \, = \, 0}^{\infty}}
    x_k^T (0)
    P_k
    x_k (0).
    \ea
    \non
    \eeq

We will show the exponential stability of the $C_o$-semigroup $T(t)$
on $L_2^n (0,1)$ by establishing convergence of the infinite sum
$\sum_{k \, = \, 0}^{\infty} x_k^T (0) P_k x_k (0)$ for each $\{ x_k
(0) \}_{k \, \in \, \bbN_0} \in l_2^n$~\cite[Lemma~5.1.2]{curzwa95}.
Let $s_m$ denote the $m$th partial sum, i.e.
    \beq
    s_m
    \; := \;
    \ds{\sum_{k \, = \, 0}^{m}}
    x_k^T (0)
    P_k
    x_k (0).
    \label{eq.sm}
    \eeq
For $l < m$ we have
    \beq
    | s_m \, - \, s_l |
    \; = \;
    \ds{\sum_{k \, = \, l \, + \, 1}^{m}}
    x_k^T (0)
    P_k
    x_k (0)
    \; \leq
    \ds{\sum_{k \, = \, l \, + \, 1}^{m}}
    \| P_k \|
    \| x_k (0) \|^2.
    \label{eq.cauchy}
    \eeq
Now, we represent $A_k$, for $k \neq 0$,  as
    \beq
    A_k
    \; = \;
    k^2
    \left(
    F_0 \; + ({1}/{k^2}) A_0
    \right),
    ~~~
    F_0
    \; := \;
    - \pi^2
    \diag \{ \obth{c_1}{\cdots}{c_n} \}
    \, < \,
    0
    \non
    \eeq
and use perturbation analysis to express $P_k$ as
    \beq
    P_k
    \; = \;
    \dfrac{1}{k^2}
    \left(
    V_0 \, + \, \dfrac{1}{k^2} V_1 \, + \, \dfrac{1}{k^4} V_2 \, +
    \, \ldots
    \right)
    \; = \;
    \dfrac{1}{k^2}
    \ds{\sum_{j \, = \, 0}^{\infty}}
    \dfrac{1}{k^{2j}} V_j
    \non
    \eeq
where
    \beq
    F_0 V_0
    \, + \,
    V_0 F_0
    \; = \;
    - \, I,
    ~~~
    F_0 V_j
    \, + \,
    V_j F_0
    \; = \;
    -
    (
    A_0^T V_{j-1}
    \, + \,
    V_{j-1} A_0
    )
    \label{eq.lyap-pert}
    \eeq
with $j \, \in \, \bbN$. Solution to~(\refeq{eq.lyap-pert}) is
determined by
    \beq
    V_0
    \; = \;
    - ({1}/{2}) F_0^{-1},
    ~~~
    V_j
    \; = \;
    \ds{\int_{0}^{\infty}}
    \mre^{F_0 t}
    (
    A_0^T V_{j-1}
    +
    V_{j-1} A_0
    )
    \mre^{F_0 t}
    \,
    \mrd t
    \non
    \eeq
which can be used to obtain
    \beq
    \ba{rcl}
    \| V_0 \|
    & \!\! = \!\! &
    1/(2 \pi^2 c_{\min})
    \\[0.1cm]
    \| V_j \|
    & \!\!  \leq \!\! &
    \| V_0 \|
    \left(
    2 \, \| A_0 \| \, \| V_0 \|
    \right)^j,
    ~~
    j \, \in \, \bbN
    \\[0.1cm]
    \| P_k \|
    & \!\! \leq \!\! &
    \dfrac{\| V_0 \|}{k^2}
    \ds{\sum_{j \, = \, 0}^{\infty}}
    \left(
    2 \, \| A_0 \| \, \| V_0 \|/k^{2}
    \right)^j.
    \ea
    \non
    \eeq
Clearly, for $k^2 > 2 \, \| A_0 \| \, \| V_0 \|$ the geometric
series in the last inequality converges. This immediately gives the
following upper bound for $\| P_k \|$:
    \beq
    \| P_k \|
    \; \leq \;
    \dfrac{\| V_0 \|}{k^2 \, - \, 2 \, \| A_0 \| \, \| V_0 \|},
    \non
    \eeq
and inequality in~(\refeq{eq.cauchy}) simplifies to
    \beq
    | s_m - s_l |
    \, \leq \,
    \dfrac{\| V_0 \|}{(l + 1)^2 \, - \, 2 \, \| A_0 \| \, \| V_0 \|}
    \ds{\sum_{k \, = \, l \, + \, 1}^{m}}
    \| x_k (0) \|^2.
    \non
    \eeq
Hence, for each $\{ x_k (0) \}_{k \, \in \, \bbN_0} \in l_2^n$
partial sum~(\refeq{eq.sm}) represents a Cauchy sequence which
guarantees convergence of \\$\sum_{k \, = \, 0}^{\infty} x_k^T (0)
P_k x_k (0)$ and consequently
    \beq
    \int_{0}^{\infty}
    \| T(t) \psi(0) \|^2
    \,
    \mrd t
    \; < \;
    \infty,
    ~~
    \forall
    \,
    \psi(0)
    \, \in \,
    {\cal D} (\ca).
    \non
    \eeq
Since ${\cal D} (\ca)$ is dense in $L_2^n (0,1)$, by an argument as
in~\cite[p.\ 51]{ban83} this inequality can be extended to all
$\psi(0) \in L_2^n (0,1)$ which implies exponential stability of
$T(t)$~\cite[Lemma~5.1.2]{curzwa95}.

    \begin{theorem}
The $C_o$-semigroup $T(t)$ generated by
operator~\refeq{eq.LRD}-\refeq{eq.DM} on $L_2^n (0,1)$ is
exponentially stable if the secant criterion~(\refeq{eq.secant}) is
satisfied.
    \label{th.heat-exp-stab}
    \end{theorem}

\subsection{The existence of a decoupled quadratic Lyapunov function}

The following theorem extends the diagonal stability result
of~\cite{arcson06} to PDEs with $r$ spatial coordinates:

    \begin{theorem}
For system~\refeq{eq.LRD}-\refeq{eq.DM} there exist a decoupled
quadratic Lyapunov function
    \beq
    V(\psi)
    \; := \;
    \inner{\psi}{D \psi}
    \; = \;
    \sum_{i \, = \, 1}^{n}
    d_i \inner{\psi_i}{\psi_i},
    ~~
    d_i \, > \, 0,
    \label{eq.Vdiag}
    \eeq
that establishes exponential stability on $L_2^n (\Omega)$ if and
only if~(\refeq{eq.secant}) holds.
    \label{th.lyap-diag}
    \end{theorem}

    \begin{proof}
We prove the theorem for a system given by
    \beq
    \psi_t
    \; = \;
    \bar{\cal A}
    \psi
    \; := \;
    C \Delta \psi
    \, + \,
    \bar{A}_0 \psi,
    \label{eq.heat-mod}
    \eeq
where $C \, := \, \diag \{ \obth{c_1}{\cdots}{c_n} \} \, > \, 0$,
and
    \beq
    \bar{A}_0
    \; := \;
    \left[
    \ba{ccccc}
    - \, 1 & 0 & \cdots & 0 & - \, \gamma_1
    \\
    \gamma_2 & - \, 1 & \ddots & & 0
    \\
    0 & \gamma_3 & - \, 1 & \ddots & \vdots
    \\
    \vdots & \ddots & \ddots & \ddots & 0
    \\
    0 & \cdots & 0 & \gamma_n & - \, 1
    \ea
    \right].
    \label{A0}
    \eeq
This is because all operators of the form~\refeq{eq.LRD} can be
obtained by acting on $\bar{A}_0$ from the left with a diagonal
matrix which does not change the existence of a decoupled quadratic
Lyapunov function. We will prove that the secant
criterion~\refeq{eq.C3} is both necessary and sufficient for the
existence of a decoupled quadratic Lyapunov function.

{\em Necessity:\/} Suppose that there exist a Lyapunov function of
the form~(\refeq{eq.Vdiag}) that establishes exponential stability
of~(\refeq{eq.heat-mod}). The derivative of~(\refeq{eq.Vdiag}) along
the solutions of~(\refeq{eq.heat-mod}) is given by
    \beq
    \ba{rcl}
    \dfrac{\mrd V(\psi)}{\mrd t}
    & \!\! = \!\! &
    \inner{\psi_t}{D \psi}
    \; + \;
    \inner{\psi}{D \psi_t}
    \; = \;
    \inner{C \Delta \psi + \bar{A}_0 \psi}{D \psi}
    \; + \;
    \inner{\psi}{D C \Delta \psi + D \bar{A}_0 \psi}
    \\[0.1cm]
    & \!\! = \!\! &
    - \, 2 \, \ds{\sum_{i \, = \, 1}^{n}} c_i d_i \inner{\nabla \psi_{i}}{\nabla \psi_{i}}
    \; + \;
    \inner{\psi}{ ( \bar{A}_0^T D + D \bar{A}_0) \psi}
    \; \leq \;
    \inner{\psi}{ ( \bar{A}_0^T D \, + \, D \bar{A}_0) \psi}
    \ea
    \non
    \eeq
where we have used Green's integral identity~\cite{eva02} with
$\psi$ satisfying the Neumann boundary conditions on $\partial
\Omega$, and the fact that $C$ and $D$ commute. The exponential
stability of~(\refeq{eq.heat-mod}) and the above expression for
${\mrd V(\psi)}/{\mrd t}$ imply that $\bar{A}_0$ is Hurwitz.
But~\refeq{eq.C3} is a necessary condition for a matrix $\bar{A}_0$
with equal diagonal entries to be Hurwitz~\cite{tysoth78}.

{\em Sufficiency:\/} Suppose that~\refeq{eq.C3} holds.
Following~\cite{arcson06} we define:
    \beq
    r
    \; := \;
    (\gamma_1 \cdots \gamma_n)^{1/n} \, > \, 0,
    ~~~
    \Gamma
    \; := \;
    \diag
    \left\{1,\ -\dfrac{\gamma_2}{r},\ \dfrac{\gamma_2\gamma_3}{r^2},\
    \cdots,\ (-1)^{n+1}\dfrac{\gamma_2\cdots \gamma_n}{r^{n-1}}
    \right\},
    ~~~
    D
    \; := \;
    \Gamma^{-2}
    \non
    \eeq
and differentiate~(\refeq{eq.Vdiag}) along the solutions
of~(\refeq{eq.heat-mod}) to obtain
    \beq
    \ba{rcl}
    \dfrac{\mrd V(\psi)}{\mrd t}
    & \!\! \leq  \!\! &
    \inner{\psi}{ ( \bar{A}_0^T D \, + \, D \bar{A}_0) \psi}
    \; =: \,
    - \inner{\psi}{ Q \psi}.
    \ea
    \non
    \eeq
If~\refeq{eq.C3} holds then $Q = Q^T$ is a positive definite
matrix~\cite{arcson06}
    \beq
    Q
    \; := \;
    - \,
    ( \bar{A}_0^T D \, + \, D \bar{A}_0)
    \; = \;
    - \,
    \Gamma^{-1}
    (
    \Gamma \bar{A}_0^T \Gamma^{-1}
    \, + \,
    \Gamma^{-1} \bar{A}_0 \Gamma
    )
    \Gamma^{-1}
    \; > \;
    0
    \non
    \eeq
and hence
    \beq
    \dfrac{\mrd V(\psi)}{\mrd t}
    \; \leq  \;
    - \lambda_{\min} (Q)
    \| \psi \|^2,
    \non
    \eeq
where $\lambda_{\min} (Q) > 0$ denotes the smallest eigenvalue of
$Q$. Upon integration, we get
    \beq
    0
    \; \leq \;
    \inner{\psi(t)}{D \psi(t)}
    \; \leq \;
    \inner{\psi(0)}{D \psi(0)}
    \, - \,
    \lambda_{\min} (Q)
    \ds{\int_{0}^{t}}
    \| \bar{T}(t) \psi(0) \|^2
    \,
    \mrd \tau
    \non
    \eeq
which yields
    \beq
    \ds{\int_{0}^{t}}
    \| \bar{T}(t) \psi(0) \|^2
    \,
    \mrd \tau
    \; \leq \;
    \dfrac{1}{\lambda_{\min} (Q)}
    \inner{\psi(0)}{D \psi(0)},
    ~~
    \forall
    \,
    t \, \geq \, 0,
    ~~
    \forall
    \,
    \psi(0)
    \, \in \,
    {\cal D} (\bar{\ca}).
    \non
    \eeq
Since ${\cal D} (\bar{\ca})$ is dense in $L_2^n (\Omega)$, the last
inequality can be extended to all $\psi(0) \in L_2^n
(\Omega)$~\cite{ban83,curzwa95}. Thus, for every $\psi(0) \in L_2^n
(\Omega)$ there is $\mu_\psi := \inner{\psi(0)}{D
\psi(0)}/\lambda_{\min} (Q)
> 0$ such that
    \beq
    \int_{0}^{\infty}
    \| \bar{T}(t) \psi(0) \|^2
    \,
    \mrd \tau
    \; \leq \;
    \mu_\psi,
    \non
    \eeq
which proves the exponential stability of
$\bar{T}(t)$~\cite[Lemma~5.1.2]{curzwa95}.
    \end{proof}

    \begin{remark}
The exponential stability of $T(t)$ in
Theorem~\ref{th.heat-exp-stab} can be also established using a
Lyapunov based approach with
    \beq
    V(\psi)
    \; = \;
    \inner{\psi}{D \psi},
    ~~
    D
    \; := \;
    \Gamma^{-2}
    \,
    \diag \{ \obth{1/a_1}{\cdots}{1/a_n} \}.
    \non
    \eeq
However, the proof of Theorem~\ref{th.heat-exp-stab} is of
independent interest because of the explicit construction of the
$C_o$-semigroup and block-diagonalization of
operator~\refeq{eq.LRD}-\refeq{eq.DM}.

    \end{remark}

\section{Extension to nonlinear reaction-diffusion equations}
    \label{NL}

We next show global asymptotic stability of the origin of the
nonlinear distributed system~\refeq{eq.RD}-\refeq{eq.DM}. This
result holds in the $L_2^n (\Omega)$ sense under the following
assumption:
    \begin{assumption}
    \label{ass.fgh}
The functions $f_i (\cdot)$, $g_i (\cdot)$, and $h_i (\cdot)$
in~\refeq{eq.RD} are continuously differentiable. Moreover,
the functions $f_i (\cdot)$ and $g_i (\cdot)$ satisfy
\refeq{eq.C1}-\refeq{eq.C3}, the functions $h_i(\cdot)$ are positive,
and the functions $g_i (\cdot)$ are nondecreasing, i.e.
    \beq
    h_i
    \; > \;
    0,
    ~~
    g_{i \sigma}
    \; := \;
    {\partial g_i}/{\partial \sigma}
    \; \geq \;
    0,
    ~~
    \forall \,
    \sigma
    \, \in \,
    \bbR.
    \tag*{(C5)}
    \label{eq.C5}
    \eeq
    \end{assumption}

A new ingredient in Assumption~\ref{ass.fgh} compared to the
properties of $f_i (\cdot)$ and $g_i (\cdot)$ in~(\refeq{bio}) is a
nondecreasing assumption on the functions $g_i (\cdot)$. This additional
assumption provides convexity of the Lyapunov function, which is
essential for establishing stability in the presence of linear
diffusion terms. For nonlinear diffusion terms we also assume that
each $h_i (\cdot)$ is a positive function.

    \begin{theorem}
    \label{thm.pde-general}
Suppose that system~\refeq{eq.RD}-\refeq{eq.DM} satisfies
Assumption~\ref{ass.fgh}. Consider the Lyapunov function candidate
    \beq
    V(\psi)
    \; = \;
    \ds{\sum_{i \, = \, 1}^{n}}
    d_i \gamma_i
    \int_{\Omega}
    \left(
    \int_0^{\psi_i(\xi)} g_i (\sigma)
    \, \mrd \sigma
    \right)
    \mrd \xi
    \non
    \eeq
where the $d_i$'s are defined as in Section~\ref{sec.linear}, and
suppose that there exists some function $\alpha (\cdot)$ of class
${\cal K}_\infty$ such that
    \beq
    V(\psi)
    \; \geq \;
    \alpha (\| \psi \|),
    ~~
    \forall \, \psi \, \in \, L_2^n(\Omega).
    \tag*{(C6)}
    \label{eq.C6}
    \eeq
Then $\psi = 0$ is a globally asymptotically stable equilibrium
point of~\refeq{eq.RD}-\refeq{eq.DM}, in the $L_2^n (\Omega)$ sense.
    \end{theorem}

    \begin{remark}[Well-posedness]
Standard arguments~(see, for example,~\cite{rot84,ama95,ama03}) can
be used to establish that~\refeq{eq.RD}-\refeq{eq.DM} has a unique
solution on $[0, \, t_{\max})$. The existence of a unique solution
on the time interval $[0, \, \infty)$ follows from the asymptotic
stability of the origin of~\refeq{eq.RD}-\refeq{eq.DM}.
    \end{remark}

    \begin{proof}
We represent the $i$th subsystem of~\refeq{eq.RD}-\refeq{eq.DM} by:
    \beq
    H_i:
    ~
    \left\{
    \ba{rcl}
    \psi_{it}
    & = &
    \nabla \cdot \left( h_i (\psi_i) \, \nabla \psi_{i} \right)
    \; - \;
    f_i(\psi_i)
    \; + \; u_i
    \\[0.1cm]
    y_i
    & = &
    g_i(\psi_i)
    \\[0.1cm]
    \psi_{i \nu}
    & = &
    0
    ~\,
    \mbox{on}
    ~\,
    \partial \Omega.
    \ea
    \right.
    \non
    \eeq
The derivative of
    \beq
    \label{pop.pde}
    V_i(\psi_i)
    \; := \;
    \gamma_i
    \int_{\Omega}
    \left(
    \int_0^{\psi_i(\xi)} g_i (\sigma)
    \, \mrd \sigma
    \right)
    \mrd \xi
    \eeq
along the solutions of $H_i$ is determined by
    \beq
    \dot{V}_i
    \; = \;
    \gamma_i
    \inner{g_i (\psi_i)}{\psi_{it}}
    \; = \;
    \gamma_i
    \inner{g_i (\psi_i)}
    { \nabla \cdot \left( h_i (\psi_i) \, \nabla \psi_{i} \right)
    \, - \,
    f_i(\psi_i)
    \, + \,
    u_i}.
    \non
    \eeq
Green's integral identity~\cite{eva02}, in combination with the
Neumann boundary conditions on $\psi_i$, can be used to obtain
    \beq
    \dot{V}_i
    \; = \;
    - \, \gamma_i \inner{g_{i \psi_i} \nabla \psi_{i}}
    {h_i \nabla \psi_{i}}
    \, - \,
    \gamma_i
    \inner{g_i}{f_i}
    \, + \,
    \gamma_i
    \inner{g_i}{u_i}.
    \non
    \eeq
Now, from~\refeq{eq.C5} we have $h_i g_{i \sigma} \geq 0$. Using
this property and the fact that $- \gamma_i f_i(\sigma) g_i(\sigma)
\leq - g_i^2 (\sigma)$ (cf.~\refeq{eq.C1}-\refeq{eq.C2}) we arrive
at
    \beq
    \dot{V}_i
    \; \leq \;
    -
    \inner{g_i}{g_i}
    \, + \,
    \gamma_i
    \inner{g_i}{u_i}
    \; = \;
    -
    \inner{y_i}{y_i}
    \, + \,
    \gamma_i
    \inner{y_i}{u_i}.
    \non
    \eeq
This upper bound on $\dot{V}_i$ and the following Lyapunov function
candidate:
    \beq
    V(\psi)
    \; := \;
    \ds{\sum_{i \, = \, 1}^{n} d_i V_i (\psi_i)}
    \non
    \eeq
yield
    \beq
    \dot{V}
    \; \leq \;
    \inner{y}{ ( \bar{A}_0^T D \, + \, D \bar{A}_0 ) y}
    \; \leq \;
    - \lambda_{\min}(Q) \| y \|^2
    \; = \;
     - \lambda_{\min}(Q)
    \ds{\sum_{i \, = \, 1}^{n}}  \| g_i \|^2.
    \label{eq.Vdot}
    \eeq
Since the $d_i$'s are defined as in Section~\ref{sec.linear}, we
have used the fact that $Q = Q^T := -(\bar{A}_0^T D \, + \, D
\bar{A}_0)$ represents a positive definite matrix (see the proof of
Theorem~\ref{th.lyap-diag}).

Now, since $V(\psi) \geq \alpha (\| \psi \|)$ for each $\psi \in
L_2^n(\Omega)$, with $\alpha (\cdot) \in {\cal K}_\infty$, for any
$\epsilon > 0$ there exist $\delta > 0$ such that $\| \psi(0) \| <
\delta$ implies $\| \psi(t) \| < \epsilon$ for all $t \geq 0$. This
follows from positive invariance of the set
    $
    \Omega_k \, := \, \{ \psi \in L_2^n(\Omega),~V(\psi) < k \},
    $
$k > 0$, and continuity of Lyapunov function $V$~\cite{hen81}.
Furthermore, $V(\psi)$ is a nonincreasing function of time bounded
below by zero and, thus, there exists a limit of $V(\psi(t))$ as
time goes to infinity. If this limit is positive then
\refeq{eq.C1},~\refeq{eq.C6}, and (\refeq{eq.Vdot}) imply the
existence of $m > 0$ such that $\sup_{t \geq 0} \dot{V} (\psi(t))
\leq - m$. But then $V(\psi(t)) \leq V(\psi(0)) - m t$ and
$V(\psi(t))$ will eventually become negative which contradicts
nonnegativity of $V(\psi(t))$, for all $t \geq 0$. Therefore, both
$V(\psi(t))$ and $\| \psi(t) \|$ converge asymptotically to zero.
>From the radial unboundedness of $V(\psi)$ (cf.\ \refeq{eq.C6}) and
the above analysis we conclude global asymptotic stability of the
origin, in the $L_2^n(\Omega)$ sense.
    \end{proof}

    \begin{remark}
The condition~\refeq{eq.C6} on $V(\psi)$ can be weakened by working
on $L_1^n(\Omega)$, in which case Jensen's inequality, applied
to~(\refeq{pop.pde}), provides the desired estimate (see
Appendix~\ref{sec.jensen}). This relaxation allows for inclusion of
many relevant nonlinearities arising in biological applications; one
such example is provided in Section~\ref{sec.example}. Using a
similar argument to the one presented in
Theorem~\ref{thm.pde-general}, the global asymptotic stability of
the origin in the $L_1^n (\Omega)$ sense can be established (with
keeping in mind that, in this case, $\inner{u}{v}$ denotes a {\em
symbol\/} for $\int_{\Omega} u^T(\xi) \, v(\xi) \, \mrd \xi$).
    \end{remark}

\section{Stability analysis for a compartmental model}
    \label{sec.compart}

An alternative to the partial differential equation
representation~\refeq{eq.RD} is a {\it compartmental model\/} which
divides the reaction into compartments that are individually
homogeneous and well-mixed,  and represents them with ordinary
differential equations. Compartmental models are preferable in
situations where reactions are separated by physical barriers such
as cell and intracellular membranes which allow limited flow between
the compartments \cite{jacquez72}. Instead of the lumped model
(\refeq{bio}) we now consider $m$ compartments where the dynamics of
the $j$th compartment, $j=2,\cdots,m-1$, are given by
    \beq
    \ba{rcl}
    \dot{x}_{j,1}
    &=&
    \mu_{j-1,1}(x_{j-1,1}-x_{j,1})
    \; - \;
    \mu_{j,1}(x_{j,1}-x_{j+1,1})
    \; - \;
    f_{1}(x_{j,1})-g_{n}(x_{j,n})
    \\[0.1cm]
    \dot{x}_{j,2}
    &=&
    \mu_{j-1,2}(x_{j-1,2}-x_{j,2})
    \; - \;
    \mu_{j,2}(x_{j,2}-x_{j+1,2})
    \; - \;
    f_{2}(x_{j,2})+g_{1}(x_{j,1})
    \\
    &\vdots&
    \\
    \dot{x}_{j,n}
    &=&
    \mu_{j-1,n}(x_{j-1,n}-x_{j,n})
    \; - \;
    \mu_{j,n}(x_{j,n}-x_{j+1,n})
    \; - \;
    f_{n}(x_{j,n})+g_{n-1}(x_{j,n-1}).
    \ea
    \tag*{(CM)}
    \label{eq.CM}
    \eeq
The functions $\mu_{j,i}(\cdot)$
$i=1,\cdots,n$, $j=1,\cdots,m-1$, represent the diffusion terms
between the compartments and possess the property
    \begin{equation}
    \sigma \mu_{j,i}(\sigma)
    \, \geq \, 0,
    ~~
    \forall \, \sigma \, \in \, \bbR.
    \tag*{(C7)}
    \label{eq.C7}
    \end{equation}
For the first and last compartments $j=1$ and $j=m$, respectively
the first and the second terms in the right-hand side
of~\refeq{eq.CM} must be dropped because $x_{0,i}$ and $x_{m+1,i}$
are not defined.

In the absence of the diffusion terms, the dynamics of the
compartments in~\refeq{eq.CM} are decoupled, and coincide with
(\refeq{bio}) which is shown in \cite{arcson06} to be globally
asymptotically stable under the conditions
\refeq{eq.C1}-\refeq{eq.C4}. The following theorem makes an
additional assumption that the function $g_i(\cdot)$ be
nondecreasing and proves that global asymptotic stability is
preserved in the presence of diffusion terms:

    \begin{theorem}
    \label{thm.compart}
Consider the compartmental model~\refeq{eq.CM}, $j=1,\dots,m$, where
for $j=1$ and $j=m$, respectively the first and the second terms in
the right-hand side of~\refeq{eq.CM} are to be interpreted as zero.
If the functions $f_{i}(\cdot)$ and $g_{i}(\cdot)$ satisfy the
conditions \refeq{eq.C1}-\refeq{eq.C4} and if, further,
$g_{i}(\cdot)$ is a nondecreasing function and $\mu_{j,i}(\cdot)$ is
as in~\refeq{eq.C7} then the origin $x_{j,i}=0$ is globally
asymptotically stable.
    \end{theorem}

\begin{proof}We first introduce the notation
    \beq
    \ba{rcl}
    x_{j}
    &  := &
    \obth{x_{j,1}}{ \cdots }{x_{j,n}}^T,
    ~~
    j \, = \, 1,\cdots,m,
    ~~
    x
    \; := \;
    \obth{x_1^T}{\cdots}{x_m^T}^T
    \\[0.1cm]
    \mu_{j}(x_j-x_{j+1})
    &  := &
    \obth{\mu_{j,1}}{\cdots}{\mu_{j,n}}^T,
    ~~
    j \, = \, 1,\cdots,m-1.
    \label{mudef}
    \ea
    \eeq
In the absence of the diffusion terms in~\refeq{eq.CM}, the reference
\cite{arcson06} constructs a Lyapunov function of the form
\begin{equation}\label{pop}
V(x_{j})=\sum_{i=1}^nd_i\gamma_i\int_0^{x_{j,i}}g_{i}(\sigma) \,
\mrd \sigma
\end{equation}
where $d_i$, $i=1,\cdots,n$, are the diagonal entries of a matrix
$D$ obtained from (\refeq{dstable}) with $A$ selected as in
(\refeq{A0}), and proves that its time derivative satisfies the
estimate
\begin{equation}\label{nodiff}
\dot{V}(x_j)\le -\epsilon \|(g_1(x_{j,1}),\cdots,g_n(x_{j,n}))\|^2
\end{equation}
for some $\epsilon>0$. In the presence of the diffusion terms
in~\refeq{eq.CM}, the estimate (\refeq{nodiff}) becomes:
    \beq
    \ba{c}
    \dot{V}(x_j)
    \; \leq \;
    - \, \epsilon
    \|(g_1(x_{j,1}),\cdots,g_n(x_{j,n}))\|^2
    \; + \;
    \dfrac{\partial V(x_j)}{\partial x_j}\mu_{j-1}(x_{j-1}-x_j)
    \; - \;
    \dfrac{\partial V(x_j)}{\partial x_j}\mu_{j}(x_{j}-x_{j+1}),
    \\[0.1cm]
    j \; = \; 2,\cdots,m-1
    \ea
    \eeq
while for $j=1$:
    \beq
    \dot{V}(x_1)
    \; \leq \;
    - \, \epsilon \|(g_1(x_{1,1}),\cdots,g_n(x_{1,n}))\|^2
    \; - \;
    \dfrac{\partial
    V(x_1)}{\partial x_1}\mu_{1}(x_{1}-x_{2})
    \label{diff1}
    \eeq
and for $j=m$:
    \beq
    \dot{V}(x_m)
    \; \le \;
    - \, \epsilon \|(g_1(x_{m,1}),\cdots,g_n(x_{m,n}))\|^2
    \; + \;
    \frac{\partial V(x_m)}{\partial x_m}\mu_{m-1}(x_{m-1}-x_m).
    \non
    \eeq
Then the Lyapunov function
\begin{equation}
\mathcal{V}(x) \; = \; \sum_{j=1}^{m}V(x_j)
\end{equation}
satisfies
    \beq
    \dot{\mathcal{V}}(x)
    \; \le \;
    - \, \epsilon \,
    \ds{\sum_{j=1}^m} \|(g_1(x_{j,1}),\cdots,g_n(x_{j,n}))\|^2
    \; - \;
    \ds{\sum_{j=1}^{m-1}} \left(\frac{\partial V(x_j)}{\partial x_j}
    \, - \,
    \frac{\partial V(x_{j+1})}{\partial x_{j+1}}
    \right)\mu_j(x_j-x_{j+1}).
    \label{re1}
    \eeq
Substituting (\refeq{mudef}) and
    \begin{equation}
    \frac{\partial V(x_j)}{\partial x_j}
    \, = \,
    \obth{d_1\gamma_1 g_1(x_{j,1})}{\cdots}
    {d_n\gamma_ng_n(x_{j,n})}
    \end{equation}
which is obtained from (\refeq{pop}), we get
    \beq
    \sum_{j=1}^{m-1}\left(\frac{\partial V(x_j)}{\partial
    x_j}-\frac{\partial V(x_{j+1})}{\partial x_{j+1}}
    \right)\mu_j(x_j-x_{j+1})
    \; = \;
    \sum_{j=1}^{m-1}\sum_{i=1}^n
    d_i\gamma_i[g_i(x_{j,i})-g_i(x_{j+1,i})]\,\mu_{j,i}(x_{j,i}-x_{j+1,i}).
    \label{re2}
    \eeq
Because $g_i(\cdot)$ is a nondecreasing function by assumption, we
note that $[g_i(x_{j,i})-g_i(x_{j+1,i})]$ possesses the same sign as
$(x_{j,i}-x_{j+1,i})$. We next recall from property~\refeq{eq.C7}
that $\mu_{j,i}(x_{j,i}-x_{j+1,i})$ also possesses the same sign as
$(x_{j,i}-x_{j+1,i})$ and, thus,
\begin{equation}
[g_i(x_{j,i})-g_i(x_{j+1,i})]\,\mu_{j,i}(x_{j,i}-x_{j+1,i})\ge 0
\end{equation}
which, according to (\refeq{re2}) and (\refeq{re1}), implies
\begin{equation}\label{re3}
 \dot{\mathcal{V}}(x)\le -\epsilon
 \sum_{j=1}^m\|(g_1(x_{j,1}),\cdots,g_n(x_{j,n}))\|^2.
\end{equation}
Because the Lyapunov function $\mathcal{V}(x)$ is proper from
property~\refeq{eq.C4} and because the right-hand side of
(\refeq{re3}) is negative definite
 from property \refeq{eq.C1}, we conclude that the origin $x=0$ is
 globally asymptotically stable. \mbox{} \hfill
\end{proof}

    \begin{remark}
Theorems \ref{thm.pde-general} and \ref{thm.compart} both rely on
the assumption that $g_i(\cdot)$ is nondecreasing, which translates
to the convexity of the Lyapunov functions (\refeq{pop.pde}) and
(\refeq{pop}). A similar convex Lyapunov function assumption has
been employed in \cite{FitHolMor97} to preserve asymptotic stability
in the presence of diffusion terms. Unlike the local result in
\cite{FitHolMor97}, however, in this paper we have established {\it
global} asymptotic stability and allowed nonlinear diffusion terms
by exploiting the specific structure of the system.
    \end{remark}

\section{An example}
    \label{sec.example}

We illustrate our main results with the analysis of a negative
feedback loop around a simple MAPK cascade model.  As described in
the introduction, MAPK cascades are functional modules, highly
conserved throughout evolution and across species, which mediate the
transmission of signals generated by receptor activation into
diverse biochemical and physiological responses involving cell cycle
regulation, gene expression, cellular metabolism, stress responses,
and other functions.  The control of MAPK and similar kinase
cascades by therapeutic intervention is being investigated as a
target for drugs, particularly in the areas of cancer and
inflammation~\cite{cohen-nrdd}. Several MAPK cascades have been
found in yeast~\cite{pg-science} and at least a dozen in mammalian
cells~\cite{karan}, and much effort is directed to the understanding
of their dynamical behavior
\cite{lauffenburger,bhalla-science,huang-ferrel}.

There are many models of MAPK cascades, with varying complexity. The
simplest class of models~\cite{heinrich02,KholodenkoNegfeedback},
using quasi-steady state approximations for enzymatic mechanisms and
a single phosphorylation site, involves a chain of three subsystems:
\begin{eqnarray*}
\label{MAPK1} \dot x_1 &=& -\,\frac{b_1 x_1}{c_1+x_1}\,+
u
\,\frac{d_1(1-x_1)}{e_1+(1-x_1)}\\
\label{MAPK2} \dot x_2 &=& -\,\frac{b_2 x_2}{c_2+x_2}\,+ \,x_1
\,\frac{d_2(1-x_2)}{e_2+(1-x_2)}\\
\label{MAPK3} \dot x_3 &=& -\,\frac{b_3 x_3}{c_3+x_3}\,+ \,x_2
\,\frac{d_3(1-x_3)}{e_3+(1-x_3)}
\end{eqnarray*}
where $u$ is an input and $x_3$ is seen as an output. The variables
$x_i$ denote the ``active'' forms of each of three proteins, and the
terms $1-x_i$ indicate the inactive forms of the respective proteins
(after nondimensionalizing and assuming that the total concentration
of each of the proteins $=1$). For example, the term $x_1
{d_2(1-x_2)}/{(e_2+(1-x_2))}$ indicates the rate at which the
inactive form of the second protein is being converted to active
form.  This rate is proportional to the concentration of the active
form of the protein $x_1$, which facilitates the conversion.
Similarly, active $x_2$ facilitates the activation of the third
protein. The first term in each of the right-hand sides models the
inactivation of the respective protein, a mechanism that proceeds at
a rate that is independent of the activation process. The saturated
form of the nonlinearities reflects the assumption that reactions
are rate limited by resources such as the amount of enzymes
available (an assumption that is not always valid). For this model,
Kholodenko proposed in \cite{KholodenkoNegfeedback} the study of
inhibitory feedback from the last to the first element,
mathematically represented by a feedback law $u = {\mu}/{(1+k
x_3)}$. See~\cite{KholodenkoNegfeedback} for a description of the
physical mechanism (an inhibitory phosphorylation of ``SOS''
protein, upstream of the system, by the last protein, p42/p44 MAPK
or ERK) that might produce this inhibition.

\begin{figure*}[ht!]
    \begin{center}
    \begin{tabular}{ccc}
    $\psi_1(\xi,t):$
    &
    $\psi_2(\xi,t):$
    &
    $\psi_3(\xi,t):$
    \\
    \includegraphics[height=1.5in,width=2.in]{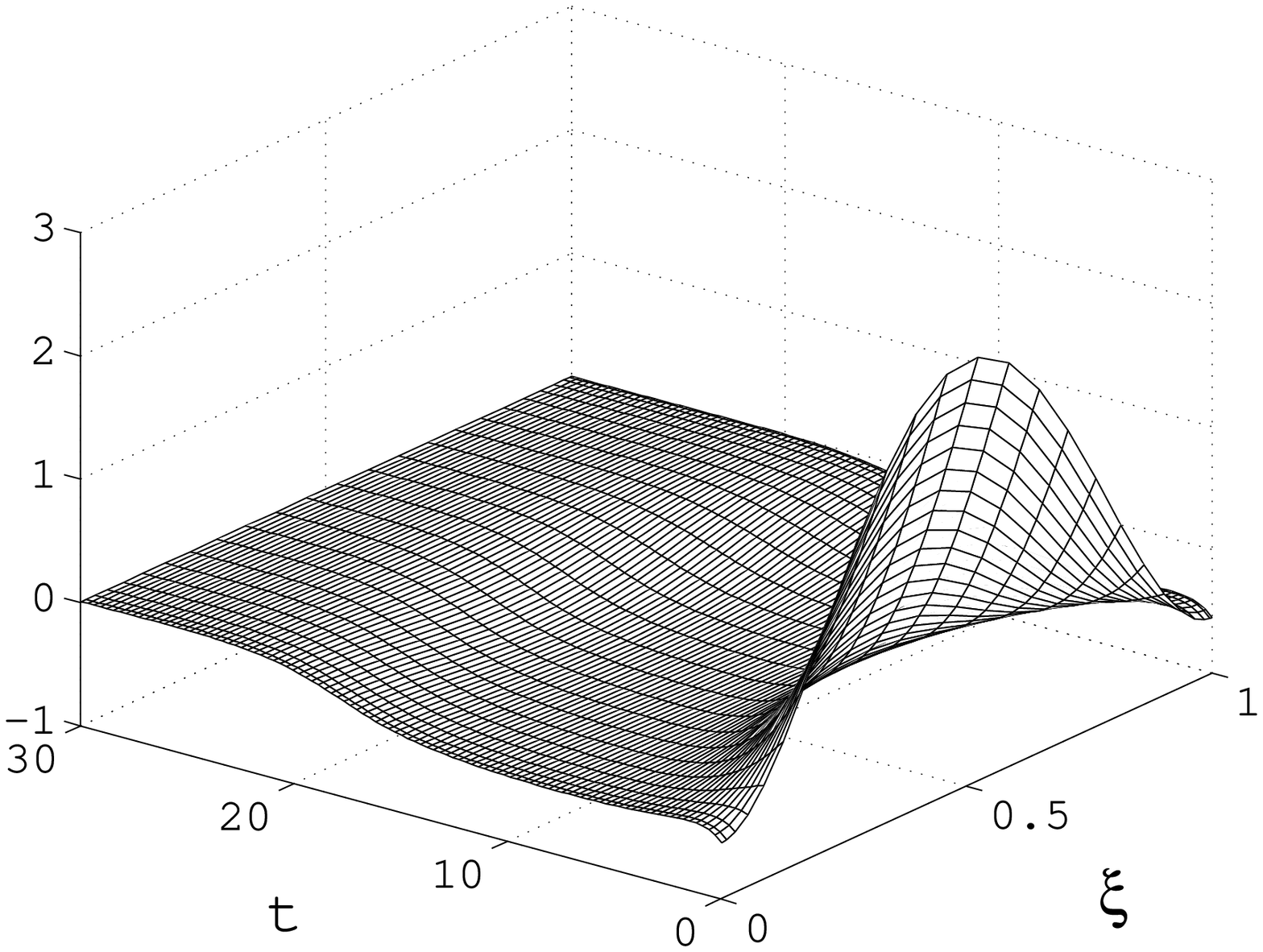}
    &
    \includegraphics[height=1.5in,width=2.in]{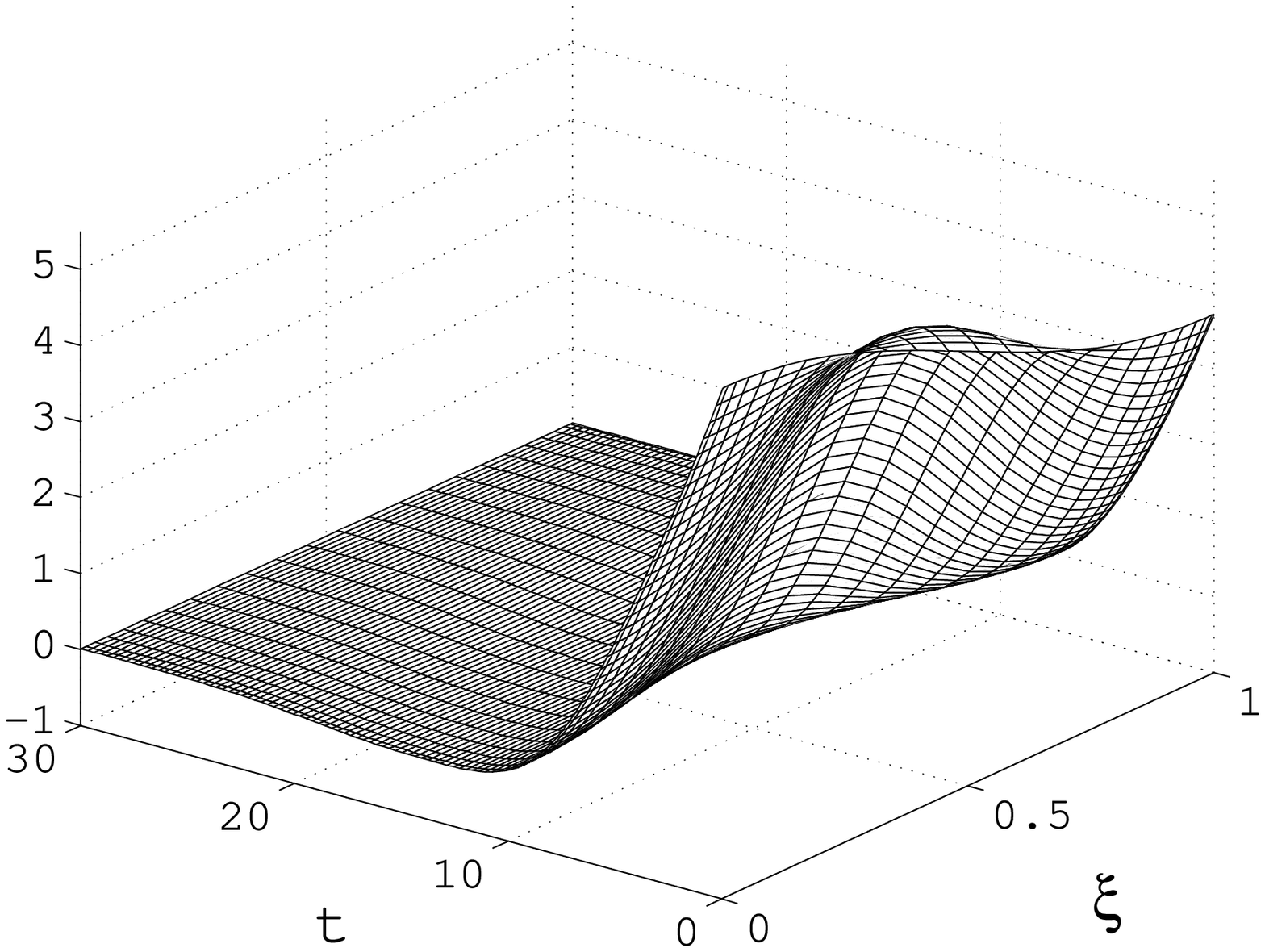}
    &
    \includegraphics[height=1.5in,width=2.in]{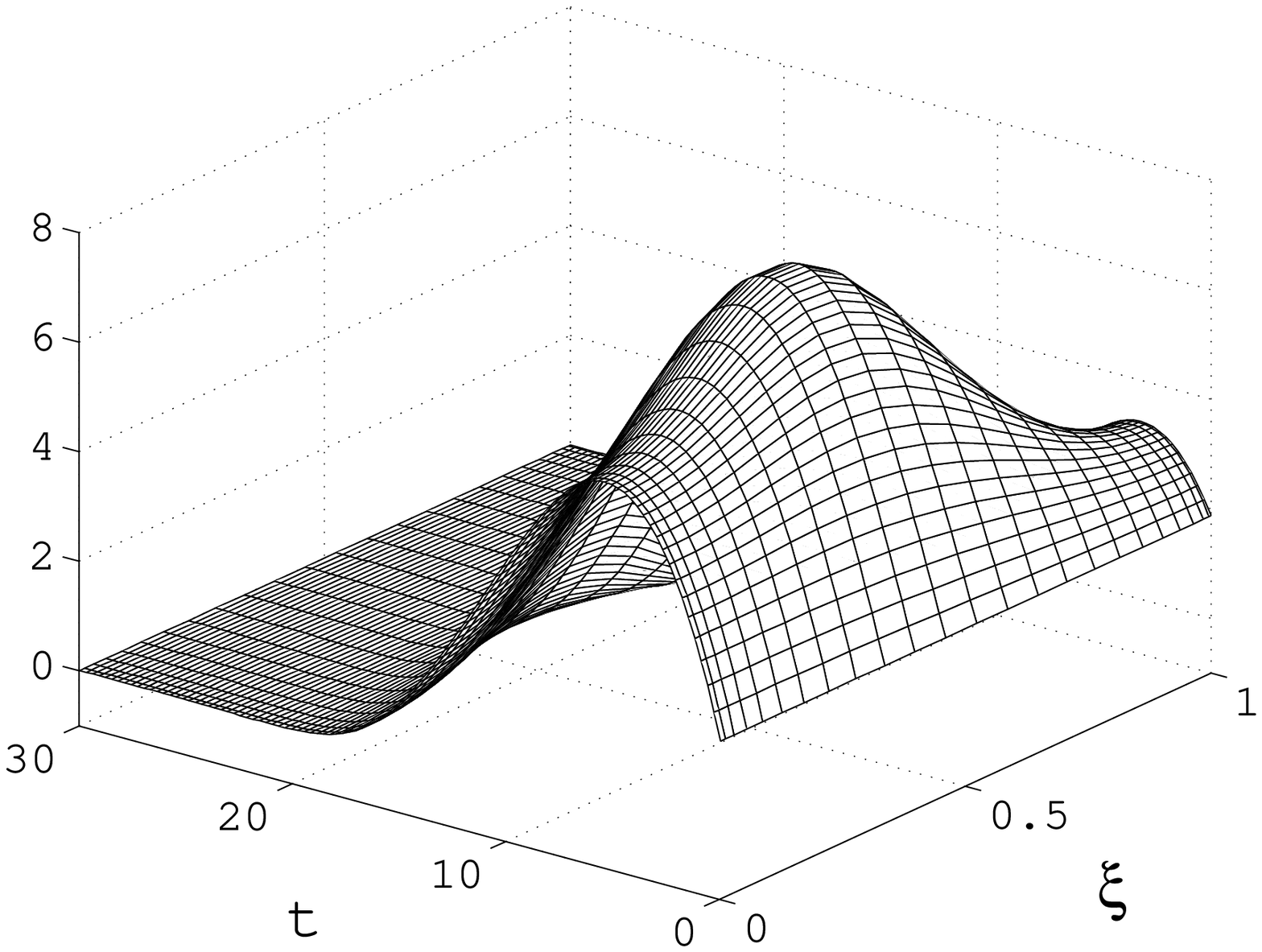}
    \end{tabular}
    \end{center}
    \caption{Plots of
    $\psi_{i} (\xi,t) \, := \, \phi_{i}(\xi,t) - \bar{\phi}$,
    $i = 1,2,3$,
    for system~\refeq{eq.EX} with
    $
    \phi(\xi,0) \, = \,
    \left[
    ~
    {16 \xi^2 (1 - \xi^2)^2}
    ~~
    {5 + \cos \, \pi \xi}
    ~~
    {2}
    ~
    \right]^T.
    $
    }
    \label{fig.mesh}
    \end{figure*}

Linearizing the system about an equilibrium, there results a linear
system to which one may apply the secant
condition~\cite{KholodenkoNegfeedback}. A linear model also arises
when considering {\it weakly activated pathways}, the behavior of
the pathway when there is only a low level of kinase
phosphorylation.  In this case, one assumes that the inactive forms
dominate: $1-x_i\approx1$; this is the analysis in~\cite{heinrich02}
and~\cite{signaling_gains04}. An intermediate case would be that in
which activations are weak but the coefficients $c_i$ are small
enough that the negative terms in the above equations cannot be
replaced by linear functions; in that case, we are lead to equations
as follows, for the closed-loop system:
    \beq
    \ba{rcl}
    \dot x_1
    &=&
    -\, \dfrac{b_1 x_1}{c_1+x_1}
    \,+ \,
    \dfrac{\mu}{1+k x_3}
    \\[0.3cm]
    \dot x_2
    &=&
    -\, \dfrac{b_2 x_2}{c_2+x_2}
    \,+ \, \,d_2x_1
    \\[0.3cm]
    \dot x_3
    &=&
    -\,\dfrac{b_3 x_3}{c_3+x_3} \,+ \, d_3 x_2
    \ea
    \label{mapk1}
    \eeq
Both linearizations, as well as this nonlinear system, can be
analyzed using our techniques. For our simulations we pick the
nonlinear model, as it is more interesting. Denoting by $\bar{x}$
the equilibrium of (\ref{mapk1}) and introducing the shifted
variable $\tilde{x}=x-\bar{x}$, we represent system (\ref{mapk1}) as
in (\ref{bio}) with
    \beq
    \ba{rcl}
    f_i(\tilde{x}_i)
    &=&
    \dfrac{b_ix_i}{c_i+x_i}
    \, - \,
    \dfrac{b_i\bar{x}_i}{c_i+\bar{x}_i},
    ~~
    i=1,2,3,
    \\[0.25cm]
    g_i(\tilde{x}_i)
    &=&
    d_{i+1}\tilde{x}_i,
    ~~
    i=1,2,
    \\[0.1cm]
    g_3(\tilde{x}_3)
    &=&
    \dfrac{\mu}{1+k \bar{x}_3}
    \, - \,
    \dfrac{\mu}{1+kx_3}.
    \ea
    \non
    \eeq
We first note that condition (C1) is satisfied because
$f_i(\tilde{x}_i)$ and $g_i(\tilde{x}_i)$, $i=1,2,3$, are strictly
increasing functions. Next, we recall from \cite[Section
6]{arcson06} that a set of gains $\gamma_i$, $i=1,2,3$, that do not
depend on the specific location of $\bar{x}$ can be obtained by
evaluating the maximum value of the slope ratio $g'_i/f_i'$ in the
interval $x_i\in [0,1]$ in which $x_i$ evolves. Upon trivial
calculations we obtain:
    \beq
    \gamma_i
    \; = \;
    \dfrac{d_{i+1}(c_i+1)^2}{b_ic_i},
    ~
    i=1,2;
    ~~
    \gamma_3
    \; = \;
    \dfrac{k\mu}{b_3c_3} \, \max \, \{c_3^2, \, \dfrac{(c_3+1)^2}{(1+k)^2}\}.
    \non
    \eeq
We pick the parameters $b_i=c_i=1$, $i=1,2,3$, $k=1$,
$d_2=d_3=\mu=0.4$, which satisfy the secant condition (C3).
(Although we chose the parameters to be $\mathcal{O}(1)$ as
in~\cite{heinrich02}, we do not claim that these are physiologically
realistic, nor are the diffusion constants that we pick below.  We
are merely interested in illustrating the theoretical results.)
Adding diffusion terms with  coefficients $h_1=h_2=h_3=0.001$, we
obtain the reaction-diffusion equations
    \beq
    \ba{rcl}
    \phi_{1t}
    &  = &
    0.001 \, \phi_{1 \xi \xi}
    \, - \,
    \dfrac{\phi_1}{1 \, + \, \phi_1}
    \, + \,
    \dfrac{0.4}{1 \, + \, \phi_3}
    \\[0.25cm]
    \phi_{2t}
    &  =  &
    0.001 \, \phi_{2 \xi \xi}
    \, - \,
    \dfrac{\phi_2}{1 \, + \, \phi_2}
    \, + \,
    0.4 \phi_1
    \\[0.25cm]
    \phi_{3t}
    &  =  &
    0.001 \, \phi_{3 \xi \xi}
    \, - \,
    \dfrac{\phi_3}{1 \, + \, \phi_3}
    \, + \,
    0.4 \phi_2
    \ea
    \tag*{(EX)}
    \label{eq.EX}
    \eeq
with the Neumann boundary conditions, $\phi_{i\xi} (0,t) =
\phi_{i\xi} (1,t) =  0$. System~\refeq{eq.EX} can be brought to the
form~\refeq{eq.RD} using the following coordinate transformation:
$\psi := \phi - \bar{\phi}$, where $\bar{\phi} \, = \,
\obth{0.5501}{0.2821}{0.1272}^T$ denotes the equilibrium point
of~\refeq{eq.EX}. Asymptotic convergence of $\phi(x,t)$ to
$\bar{\phi}$ is illustrated in Fig.~\ref{fig.mesh}. A spatial
discretization of the diffusion operator with Neumann boundary
conditions is obtained using a Matlab Differentiation Matrix
Suite~\cite{weired00}.

\section{Concluding remarks}

We identify a class of systems with a cyclic interconnection
structure in which addition of diffusion does not have a
destabilizing effect. For these systems, we demonstrate global
stability if the ``secant'' criterion is satisfied. In the linear
case, we show that the secant condition is necessary and sufficient
for the existence of a decoupled Lyapunov function, which extends
the diagonal stability result~\cite{arcson06} to spatially
distributed systems. For reaction-diffusion equations with
nondecreasing coupling nonlinearities, we establish global asymptotic
stability of the origin. Under some fairly mild assumptions, we also
allow for nonlinear diffusion terms by exploiting the specific
structure of the system.

    \appendix

\subsection{Relaxation of condition~\refeq{eq.C6}}
    \label{sec.jensen}

Let us represent $V_i(\psi_i)$ in~(\ref{pop.pde}) by
    \beq
    V_i(\psi_i)
    \; := \;
    \gamma_i
    \ds{\int_{\Omega}}
    p_i (\psi_i(\xi))
    \,
    \mrd \xi,
    ~~~
    p_i (s)
    \; := \;
    \ds{\int_0^{s} g_i (\sigma)}
    \, \mrd \sigma.
    \non
    \eeq
and let $\Omega_p$ (respectively, $\Omega_m$) denote the set of
points in $\Omega$ where $\psi_i(\xi)$ is positive (respectively,
negative), i.e.
    \beq
    \Omega_p
    \; := \;
    \left\{
    \xi \, \in \, \Omega,
    ~
    \psi_i(\xi) \, > \, 0
    \right\},
    ~~~
    \Omega_m
    \; := \;
    \left\{
    \xi \, \in \, \Omega,
    ~
    \psi_i(\xi) \, < \, 0
    \right\}.
    \non
    \eeq
Then, $V_i(\psi_i)$ can be rewritten as
    \beq
    V_i(\psi_i)
    \; := \;
    \gamma_i
    \ds{\int_{\Omega_p}}
    p_{ip} (| \psi_i(\xi) |)
    \,
    \mrd \xi
    \; + \;
    \gamma_i
    \ds{\int_{\Omega_m}}
    p_{im} (| \psi_i(\xi) |)
    \,
    \mrd \xi
    \non
    \eeq
where
    \beq
    p_{ip} (s)
    \; := \;
    \ds{\int_0^{s} g_i (\sigma)}
    \, \mrd \sigma,
    ~
    p_{im} (s)
    \; := \;
    \ds{\int_0^{-s} g_i (\sigma)}
    \, \mrd \sigma,
    ~
    s \, > \, 0.
    \non
    \eeq
We observe that the first two derivatives of the functions $p_{ip}$ and
$p_{im}$, respectively, satisfy
    \beq
    \ba{c}
    \{
    p'_{ip} (s)
    \; = \;
    g_i (s)
    \; > \; 0,
    ~~
    p''_{ip} (s)
    \; = \;
    g'_i (s)
    \; \geq \; 0
    \}
    \\[0.15cm]
    \{
    p'_{im} (s)
    \; = \;
    - g_i (-s)
    \; > \; 0,
    ~~
    p''_{im} (s)
    \; = \;
    g'_i (-s)
    \; \geq \; 0
    \}
    \ea
    \non
    \eeq
which implies that both these functions are of class ${\cal
K}_\infty$ and convex.  Using convexity, we may apply Jensen's
inequality~\cite{eva02} to obtain
    \beq
    V_i(\psi_i)
    \, \geq \,
    \gamma_i
    \left(
    |\Omega_p| \, p_{ip} ( \| \psi_i \|_1/|\Omega_p| )
    \, + \,
    |\Omega_m| \, p_{im} ( \| \psi_i \|_1/|\Omega_m| )
    \right)
    \non
    \eeq
where $|\Omega_r|$ denotes the measure of set $\Omega_r$, and $\|
\psi_i \|_1$ is the $L_1(0,1)$-norm of $\psi_i$. Since $p_{ip}$ and
$p_{im}$ are ${\cal K}_\infty$-functions we conclude that
condition~\refeq{eq.C6} on $V(\psi)$ always holds if the underlying
state-space is $L_1^n(\Omega)$ (that is, there exists some function
$\alpha (\cdot)$ of class ${\cal K}_\infty$ such that
    $
    V(\psi)
    \, \geq \,
    \alpha (\| \psi \|_1),
    $
    $
    \forall \, \psi \, \in \, L_1^n(\Omega)
    ).
    $

\end{document}